# THE SHATTERING DIMENSION OF SETS OF LINEAR FUNCTIONALS

By Shahar Mendelson[1] and Gideon Schechtman[2]

*Australian National University and Weizmann Institute*

We evaluate the shattering dimension of various classes of linear functionals on various symmetric convex sets. The proofs here relay mostly on methods from the local theory of normed spaces and include volume estimates, factorization techniques and tail estimates of norms, viewed as random variables on Euclidean spheres. The estimates of shattering dimensions can be applied to obtain error bounds for certain classes of functions, a fact which was the original motivation of this study. Although this can probably be done in a more traditional manner, we also use the approach presented here to determine whether several classes of linear functionals satisfy the uniform law of large numbers and the uniform central limit theorem.

**1. Introduction.** Combinatorial dimensions, such as the Vapnik–Chervonenkis dimension, and the shattering dimension, are parameters which measure the richness of a given class of functions. The Vapnik–Chervonenkis dimension (VC dimension) of a class of $\{0, 1\}$-valued functions is the largest dimension of a combinatorial cube that can be found in a coordinate projection of the class, that is, in a restriction of the class to a finite subset of the domain. In this article we focus on a real valued analog of the VC dimension, called the shattering dimension; it is a scale sensitive parameter that measures the largest dimension of a "cube" of a given side length that can be found in a coordinate projection of the class.

DEFINITION 1.1. For every $\varepsilon > 0$, a set $\sigma = \{x_1, \ldots, x_n\} \subset \Omega$ is said to be $\varepsilon$-shattered by a set $F$ of functions on $\Omega$ if there is some function $s : \sigma \to \mathbb{R}$, such that for every $I \subset \{1, \ldots, n\}$, there is some $f_I \in F$ for which

Received October 2002; revised June 2003.
[1] Supported by the ARC.
[2] Supported by the ISF.

$f_I(x_i) \geq s(x_i) + \varepsilon$ if $i \in I$, and $f_I(x_i) \leq s(x_i) - \varepsilon$ if $i \notin I$. The shattering dimension of $F$ is the function

$$\mathrm{VC}(\varepsilon, F, \Omega) = \sup\{|\sigma| \, \sigma \subset \Omega, \ \sigma \text{ is } \varepsilon\text{-shattered by } F\}.$$

$f_I$ is called the shattering function of the set $I$ and the set $\{s(x_i) | x_i \in \sigma\}$ is called a witness to the $\varepsilon$-shattering. In cases where the underlying space is clear we denote the shattering dimension by $\mathrm{VC}(\varepsilon, F)$.

In this article we evaluate the shattering dimension of various classes of linear functionals on various symmetric convex sets. Before describing the actual results obtained, we would like to describe the way one applies such estimates to obtain results concerning the uniform law of large numbers and the uniform central limit theorem (CLT), as well as error bounds in statistical learning theory.

Combinatorial dimensions have been frequently used in the theory of empirical processes, mostly in the context of the uniform law of large numbers and the uniform CLT. Recall the definition of the uniform law of large numbers, also known as the uniform Glivenko–Cantelli condition.

Definition 1.2.   Let $F$ be a class of functions. We say that $F$ is a uniform Glivenko Cantelli class (uGC class) if for every $\varepsilon > 0$,

$$(1.1) \qquad \lim_{n \to \infty} \sup_\mu \mathrm{Pr}\left\{ \sup_{f \in F} \left| \mathbb{E}_\mu f - \frac{1}{n} \sum_{i=1}^{n} f(X_i) \right| \geq \varepsilon \right\} = 0,$$

where $(X_i)_{i=1}^{\infty}$ are independent random variables distributed according to $\mu$.

Let us remark that in this article we ignore the question of measurability, since only mild assumptions on the class, such as admissibility, are required to resolve this issue (see [6] for further details). Moreover, in all the cases we explore, it suffices to consider the supremum over a countable dense set, and, thus, the measurability issue does not arise.

Vapnik and Chervonenkis proved that (under mild measurability assumptions) a class of binary-value functions is a uGC class if and only if it has a finite VC dimension [23], and this result was extended in [1] to the real-valued case, where it was shown that a class of uniformly bounded functions is a uGC class if and only if $\mathrm{VC}(\varepsilon, F)$ is finite for every $\varepsilon > 0$ (see also [7] for a related earlier characterization of uGC classes of functions).

The shattering dimension can be used to obtain the tail bounds needed in (1.1), using the following line of argumentation. The starting point is a version of Talagrand's inequality (originally proved in [20]) due to Bousquet.

**Theorem 1.3 ([5]).** *Let $F$ be a class of functions defined on a probability space $(\Omega, \mu)$ such that $\sup_{f \in F} \|f\|_\infty \leq 1$. Let $(X_i)_{i=1}^n$ be independent random variables distributed according to $\mu$, put $\sigma^2 \geq \sup_{f \in F} \operatorname{Var}[f(X_1)]$ and set $Z = \sup_{f \in F} |\sum_{i=1}^n (f(X_i) - \mathbb{E}_\mu f)|$. Then, for every $x > 0$,*

$$\Pr\{Z \geq \mathbb{E}Z + x\} \leq \exp\left(-vh\left(\frac{x}{v}\right)\right), \tag{1.2}$$

*where $v = n\sigma^2 + 2\mathbb{E}Z$ and $h(x) = (1 + x)\log(1 + x) - x$. Moreover, for every $x > 0$,*

$$\Pr\left\{Z \geq \mathbb{E}Z + \sqrt{2xv} + \frac{x}{3}\right\} \leq e^{-x}. \tag{1.3}$$

In order to apply this result and obtain uniform deviation estimates one needs to bound $\mathbb{E}Z$. By symmetrization,

$$\mathbb{E}_\mu Z \leq 2\mathbb{E}_{\mu \times \varepsilon} \sup_{f \in F} \left| \sum_{i=1}^n \varepsilon_i f(X_i) \right|, \tag{1.4}$$

where $(\varepsilon_i)_{i=1}^n$ are independent Rademacher random variables (i.e., take the values $\pm 1$ with probability $1/2$ each). It turns out that the Rademacher averages on the right-hand side of (1.4) can be estimated in terms of the *empirical covering numbers*.

If $(Y, d)$ is a metric space and $F \subset Y$, then for every $\varepsilon > 0$, $N(\varepsilon, F, d)$ denotes the minimal number of open balls (with respect to the metric $d$) needed to cover $F$.

**Definition 1.4.** For every class $F$ let the empirical covering numbers be

$$N(\varepsilon, F, n) = \sup_{\mu_n} N(\varepsilon, F, L_2(\mu_n)),$$

where the supremum is taken with respect to all empirical measures $n^{-1} \sum_{i=1}^n \delta_{x_i}$ supported on $n$ points. $\log N(\varepsilon, F) = \sup_n \log N(\varepsilon, F, n)$ is called the uniform $L_2$ entropy of $F$.

The following result shows that the uniform entropy can be bounded via the combinatorial parameters.

**Theorem 1.5 ([16]).** *There are absolute constants $K$ and $c$ such that for any class $F$ which consists of functions bounded by 1 and every $0 < \varepsilon < 1$,*

$$N(\varepsilon, F) \leq \left(\frac{2}{\varepsilon}\right)^{K \cdot \operatorname{VC}(c\varepsilon, F)}.$$

Combining Theorem 1.5 with a chaining argument, one can bound the Rademacher averages of (1.4) and, thus, $\mathbb{E}Z$ and obtain the necessary deviation estimates. In all the examples presented in the sequel we will establish upper bounds on the shattering dimension which are polynomial in $1/\varepsilon$, and in that case, the following holds.

THEOREM 1.6 ([15]).   *Let $F$ be a class of functions bounded by 1, and set $Z$ to be as in Theorem 1.3. Assume that there are $\gamma \geq 1$ and $0 < p < \infty$ such that $\mathrm{VC}(\varepsilon, F) \leq \gamma \varepsilon^{-p}$. Then*

$$\mathbb{E}Z \leq C_p \gamma^{1/2} \begin{cases} \sqrt{n}, & \text{if } 0 < p < 2, \\ \sqrt{n}\log^{3/2} n, & \text{if } p = 2, \\ n^{1-1/p}\log^{1/p} n, & \text{if } p > 2, \end{cases}$$

*where $C_p$ are constants which depend only on $p$.*

We now turn to the description of the connection between bounds on the shattering dimension and error bounds used in the analysis of regression problems in nonparametric statistics and, more recently, in Learning Theory. In both applications, combinatorial parameters have played an important role. In the context of Learning Theory, they were used to estimate the size of a random sample needed to construct an almost optimal approximation of an unknown target function by an element in a fixed class of functions, where the given data are a sample $(X_i)_{i=1}^n$ and the values of the target on the sample [2, 15]. Such an error bound which is based on the shattering dimension is presented in the next theorem, which was adapted from [4].

THEOREM 1.7.   *Let $(\Omega, \mu)$ be a probability space, let $F$ be a class of measurable functions on $\Omega$ with ranges in $[-1, 1]$ and assume that there is a constant $B \geq 1$ such that for every $f \in F$, $\mathbb{E}_\mu f^2 \leq B\mathbb{E}_\mu f$. If $(X_i)_{i=1}^n$ are independent random variables distributed according to $\mu$, then for any $x > 0$, there is a set of probability larger than $1 - 2e^{-x}$, on which for any $f \in F$,*

$$\mathbb{E}_\mu f \leq \frac{2}{n}\sum_{i=1}^n f(X_i) + C\left(\frac{I}{\sqrt{n}} + \frac{Bx}{n}\right),$$

*where $C$ is an absolute constant and*

$$(1.5) \qquad I = \int_0^1 \sqrt{\mathrm{VC}(\varepsilon, F, \{X_1, \ldots, X_n\})\log\left(\frac{1}{\varepsilon}\right)}\, d\varepsilon.$$

Let us mention that it is possible to obtain error bounds even in some cases when $I = \infty$ [14], and that in [4], error bounds with faster rates of convergence than $1/\sqrt{n}$ were established in the same setup.

The analysis of the shattering dimension of classes of linear functionals we present is based on methods from the local theory of normed spaces. We show that for such classes the shattering dimension is determined by the geometry of the class and the domain, which is expressed by the ability to factor a certain operator through $\ell_1^n$. First, we investigate in Section 3 the case when $\Omega$ is the unit ball of some Banach space and $F$ is the dual unit ball. We show that if $X$ is infinite dimensional and $B_X$ is the unit ball of $X$, the shattering dimension $\mathrm{VC}(\varepsilon, B_{X^*}, B_X)$ is determined by the *Rademacher type* of $X$. In Section 4, which contain the main new results of this article, we use a volumetric argument and establish estimates on the shattering dimension when both the class and the domain are finite-dimensional convex and symmetric sets. We then compute the shattering dimension of the unit ball in $\ell_q^n$ when considered as functions on the unit ball of $\ell_p^n$, $1 \leq p, q \leq \infty$, and show that in many cases the volumetric approach yields sharp bounds. For example, we prove in Theorem 4.12 that for $1 \leq p < q \leq \infty$ and for $F = B_{q'}$, the unit ball of $\ell_{q'}$ ($q'$ is the conjugate index to $q$), and $\Omega = B_p$, $\mathrm{VC}(\varepsilon, F, \Omega)$ is given, up to constants depending only on $p$ and $q$, by the following expressions:

$$\mathrm{VC}(\varepsilon, F, \Omega) \sim_{p,q} \begin{cases} \varepsilon^{-q/(q-1)}, & \text{if } 1 \leq p \leq 2, \\ \varepsilon^{-1/(1/2+1/p-1/q)}, & \text{if } 2 < p \leq \infty. \end{cases}$$

Section 5 is devoted to computation of the shattering dimension of the image of the unit ball of $\ell_1^n$ under a linear transformation.

The applications of the estimates of the shattering dimensions to the determination of whether some classes of functionals satisfy the uniform law of large numbers or the uniform CLT are scattered through Sections 3 and 4. In Section 3 we give, among other things, a new proof for a result from [7] giving a necessary and sufficient condition for the unit ball $B_{X^*}$ of a dual Banach space to be uGC class on $B_X$. In Section 4 our results are used to investigate the following problem: consider the unit ball of $\ell_q$, denoted by $B_q$, as functions on the unit ball of $\ell_p$. Does this class of functions satisfy the uniform CLT on this domain? In general, one can show that for any infinite-dimensional Banach space $X$, $F = B_{X^*}$, does not satisfy the uniform CLT on the domain $\Omega = B_X$. Although this can probably be deduced from earlier contributions, we show, as an application of the methods presented here, that whenever $p < q$, $F = B_{q'}$ satisfies the uniform CLT on the domain $\Omega = B_p$.

## 2. Preliminaries.

Throughout all absolute constants are denoted by $c$, $C$ or $K$. Their values may change from line to line or even within the same line. $c_\varphi$, $C_\varphi$ denote constants which depend only on the parameter $\varphi$ (which is usually a real number $p$ or a couple of real numbers $p, q$), and $a \sim_\varphi b$ means that $c_\varphi b \le a \le C_\varphi b$. If the constants are absolute, we use the notation $a \sim b$. Given a real Banach space $X$, let $B_X$ or $B(X)$ be the unit ball of $X$. The dual of $X$, denoted by $X^*$, consists of all the bounded linear functionals on $X$, endowed with the norm $\|x^*\| = \sup_{\|x\|=1} |x^*(x)|$. For every integer $n$, we fix the Euclidean structure $\langle \cdot, \cdot \rangle$ on $\mathbb{R}^n$ with an orthonormal basis denoted by $(e_i)_{i=1}^n$.

A set $K$ is called symmetric if the fact that $x \in K$ implies that $-x \in K$. The symmetric convex hull of $K$, denoted by $\mathrm{absconv}(K)$, is the convex hull of $K \cup -K$.

If $K \subset \mathbb{R}^n$ is bounded, convex and symmetric with a nonempty interior, then $K$ is a unit ball of a norm denoted by $\| \cdot \|_K$. It is possible to show that the *polar* of $K$, defined by

$$K^o = \left\{ x \in \mathbb{R}^n \,\Big|\, \sup_{k \in K} \langle k, x \rangle \le 1 \right\},$$

is the unit ball of the dual space of $(\mathbb{R}^n, \| \cdot \|_K)$. In the sequel we shall abuse notation and denote by $K$ the normed space whose unit ball is $K$. From here on, a ball will be a bounded, convex and symmetric subset of $\mathbb{R}^n$, with a nonempty interior.

If $1 \le p < \infty$, let $\ell_p^n$ be $\mathbb{R}^n$ endowed with the norm $\|\sum_{i=1}^n a_i e_i\|_p = (\sum_{i=1}^n |a_i|^p)^{1/p}$. $\ell_\infty^n$ is $\mathbb{R}^n$ endowed with the norm $\|\sum_{i=1}^n a_i e_i\|_\infty = \sup_i |a_i|$. $B_p^n$ is the unit ball of $\ell_p^n$, and for every $1 \le p \le \infty$, $(B_p^n)^o = B_{p'}^n$, where $1/p + 1/p' = 1$. In this case, $p'$ is called the conjugate index of $p$.

2.1. *Volume estimates.* As stated above, we can identify $\ell_2^n$ with $\mathbb{R}^n$. Hence, $\ell_2^n$ is endowed with the $n$-dimensional Lebesgue measure, denoted by $|\cdot|$. Let $GL_n$ be the set of invertible operators $T : \mathbb{R}^n \to \mathbb{R}^n$, and note that for every measurable set $A \subset \mathbb{R}^n$ and every $T \in GL_n$, $|TA| = |\det(T)| \|A|$. We say that a set $A \subset \mathbb{R}^n$ is an ellipsoid if there is some $T \in GL_n$, such that $A = T B_2^n$.

It will be useful to determine the volume of the balls $B_p^n$ and the volume of their sections. First, let us mention the following well-known fact.

THEOREM 2.1 ([19]).   *There are absolute constants $C$ and $c$ such that for every integer $n$ and every $1 \le p \le \infty$,*

$$cn^{-1/p} \le |B_p^n|^{1/n} \le Cn^{-1/p}.$$

Unlike the clear structure of sections of $B_2^n$, the geometry of sections of $B_p^n$ is far less obvious. The following result, due to Meyer and Pajor [17], bounds the volume of $k$-dimensional sections of $B_p^n$.

THEOREM 2.2. *For every $k$-dimensional subspace $E \subset \mathbb{R}^n$ and every $1 \leq p \leq q \leq \infty$,*

$$\frac{|B_p^n \cap E|}{|B_p^k|} \leq \frac{|B_q^n \cap E|}{|B_q^k|}.$$

By selecting $q = 2$, it follows that for $1 \leq p \leq 2$, the volume of any $k$-dimensional section of $B_p^n$ is smaller than the volume of $B_p^k$. Similarly, by taking $p = 2$, the volume of any $k$-dimensional section of $B_q^n$ for $2 \leq q \leq \infty$ is larger than the volume of $B_q^k$.

REMARK 2.3. A similar result holds in the infinite-dimensional case. In particular, it follows that for any $1 \leq p \leq q \leq \infty$ and for any $n$-dimensional subspace $E$,

$$(2.1) \qquad \left( \frac{|B_p \cap E|}{|B_q \cap E|} \right)^{1/n} \leq C_{p,q} n^{1/q - 1/p}.$$

An important fact about the volume of balls are the Santaló and inverse Santaló inequalities.

THEOREM 2.4. *There is an absolute constant $c$ such that for every integer $n$ and every ball $K \subset \mathbb{R}^n$,*

$$c \leq \left( \frac{|K||K^o|}{|B_2^n|^2} \right)^{1/n} \leq 1.$$

The upper bound was established by Santaló, while the lower bound is due to Bourgain and Milman. The proof of both results can be found in [19].

One of the tools used in modern convex geometry is the notion of volume ratios. The idea is to compare the volume of a given ball with the "best" possible volume of an ellipsoid contained in it, since this may be used to understand "how close" the norm induced by the ball is to a Euclidean structure.

DEFINITION 2.5. For every ball $K \subset \mathbb{R}^n$, the volume ratio of $K$ is

$$\mathrm{vr}(K) = \inf \left( \frac{|K|}{|T B_2^n|} \right)^{1/n},$$

where the infimum is taken with respect to all $T \in GL_n$ such that $T B_2^n \subset K$.

The external volume ratio is defined as

$$\mathrm{evr}(K) = \inf \left( \frac{|T B_2^n|}{|K|} \right)^{1/n},$$

where the infimum is with respect to all $T \in GL_n$ such that $K \subset T B_2^n$.

It is possible to show [19] that both infimums in the definition above are uniquely attained. Hence, for every ball $K \subset \mathbb{R}^n$, there is an ellipsoid of maximal volume contained in $K$ and an ellipsoid of minimal volume containing $K$. The ellipsoid of maximal volume contained in $K$ is denoted by $\mathcal{E}_K$, and the ellipsoid of minimal volume containing $K$ is denoted by $\tilde{\mathcal{E}}_K$. Note that for every ball $K$, $\mathcal{E}_K^o = \tilde{\mathcal{E}}_{K^o}$.

It follows from the definitions that if $K$ is an ellipsoid, then $\mathrm{vr}(K) = \mathrm{evr}(K) = 1$. Moreover, it is known that for every ball $K \subset \mathbb{R}^n$, $\mathrm{vr}(K) \leq \sqrt{n}$. More precisely, the volume ratio of $\ell_\infty^n$, which is of the order of $\sqrt{n}$, is the worst possible.

**THEOREM 2.6 ([3]).**  *For every integer $n$,*
$$\mathrm{vr}(K) \leq \mathrm{vr}(B_\infty^n) = \frac{4}{|B_2^n|^{1/n}}.$$

Another result we require is an estimate on the volume ratios of projections of $\ell_p$.

**THEOREM 2.7 ([12]).**  *For every integer $n$,*
$$\sup_{E \subset \mathbb{R}^n} \mathrm{vr}(P_E B_p) \underset{p}{\sim} \begin{cases} 1, & 1 < p \leq 2, \\ n^{1/2 - 1/p}, & 2 < p \leq \infty, \end{cases}$$
*where the supremum is taken with respect to all the projections onto $n$ dimensional subspaces of $\ell_p$.*

A different notion of volume ratios is the *cubic ratios* which was introduced by Ball [3]. For every ball $K \subset \mathbb{R}^n$, let
$$\mathrm{cr}(K) = \inf_{T \in GL_n, K \subset TB_\infty^n} \left( \frac{|TB_\infty^n|}{|K|} \right)^{1/n}.$$

**LEMMA 2.8 ([3]).**  *There are absolute constants $c$ and $C$ such that for every integer $n$ and every ball $K \subset \mathbb{R}^n$,*
$$c\sqrt{n} \leq \mathrm{vr}(K)\mathrm{cr}(K) \leq C\sqrt{n}.$$

Finally, we can define the volume numbers of an operator. We follow the definition used by Gordon and Junge [11, 12].

**DEFINITION 2.9.**  Given Banach spaces $X$ and $Y$, an operator $T : X \to Y$ and an integer $n$, let the $n$th volume number of $T$ be
$$v_n(T) = \sup \left\{ \left( \frac{|T(B_X \cap E)|}{|B_Y \cap F|} \right)^{1/n} \middle| E \subset X, \ T(E) \subset F \subset Y, \ \dim E = \dim F = n \right\}.$$

Note that if $T$ is of rank smaller than $n$, $v_n(T) = 0$. Also, it is clear that the volume numbers are submultiplicative, that is, $v_n(T_1 T_2) \leq v_n(T_1) v_n(T_2)$. If $T$ is an operator between Hilbert spaces, then the volume numbers may be calculated using the eigenvalues $(\lambda_i)$ of $\sqrt{T^* T}$ (which are arranged in a nonincreasing order). In that case, for every integer $n$, $v_n(T) = (\prod_{i=1}^n \lambda_i)^{1/n}$.

Another example in which the volume numbers may be estimated is for the formal identity operator $id : \ell_p^m \to \ell_q^m$. By Theorem 2.2 it is evident that for every $n \leq m$ and any $1 \leq p \leq q \leq \infty$,

$$(2.2) \quad v_n(id_{\ell_p^m \to \ell_q^m}) = \sup_{\dim E = n} \frac{|B_p^m \cap E|^{1/n}}{|B_q^m \cap E|^{1/n}} \leq \frac{|B_p^n|^{1/n}}{|B_q^n|^{1/n}} \leq C_{p,q} n^{1/q - 1/p}$$

and clearly also

$$(2.3) \quad v_n(id_{\ell_p^m \to \ell_q^m}) \geq \frac{|B_p^n|^{1/n}}{|B_q^n|^{1/n}}.$$

In general, if $p \geq q$, then

$$(2.4) \quad v_n(id_{\ell_p^m \to \ell_q^m}) \leq \|id\|_{\ell_p^m \to \ell_q^m} = m^{1/q - 1/p},$$

and this estimate is optimal, at least in cases where $n$ divides $m$. To see this, let $k = m/n$ and for $j = 1, \ldots, n$, let $v_j = \sum_{i=1}^k e_{j + k(i-1)}$. Note that for each $r$, $\mathrm{span}\{v_1, \ldots, v_n\} \cap B_r^m = E \cap B_r^m$ has volume $(m/n)^{1/2 - 1/r} |B_r^n|$. Thus,

$$\frac{|B_p^m \cap E|^{1/n}}{|B_q^m \cap E|^{1/n}} = \left(\frac{m}{n}\right)^{1/q - 1/p} \cdot \frac{|B_p^n|^{1/n}}{|B_q^n|^{1/n}} = C_{p,q} m^{1/q - 1/p},$$

proving that the bound on the volume numbers is tight.

2.2. *The uniform CLT.* The fact that the shattering dimension can be used to bound the uniform entropy will enable us to show that some classes of functionals satisfy the *uniform CLT*. Recall that a sequence of measures $\nu_n$ converges to $\nu$ in law in $\ell_\infty(F)$ if for every bounded and continuous function $H : \ell_\infty(F) \to \mathbb{R}$, $\mathbb{E}^* H(\nu_n) \to \mathbb{E}^* H(\nu)$, where $\mathbb{E}^*$ denotes the outer expectation.

DEFINITION 2.10 [6]. Let $F \subset B(L_\infty(\Omega))$, set $P$ to be a probability measure on $\Omega$ and assume $G_P$ to be a Gaussian process indexed by $F$, which has mean 0 and covariance

$$\mathbb{E} G_P(f) G_P(g) = \int f g \, dP - \int f \, dP \int g \, dP.$$

$F$ is called a universal Donsker class if for any probability measure $P$, the law $G_P$ is tight in $\ell_\infty(F)$ and $\nu_n^P = n^{1/2}(P_n - P) \in \ell_\infty(F)$ converges in law to $G_P$ in $\ell_\infty(F)$, where $P_n$ is a random empirical measure selected according to $P$, that is, $P_n = \frac{1}{n}\sum_{i=1}^{n}\delta_{X_i}$, where $(X_i)_{i=1}^{n}$ are independent random variables distributed according to $P$.

Stronger than the universal Donsker property is the uniform Donsker property, which is the uniform version of the CLT. For such classes, $\nu_n^P$ *converges* to $G_P$ uniformly in $P$ in some sense (see [6, 22] for more details). The following result of Giné and Zinn [8] is a relatively simple characterization of uniform Donsker classes.

For every probability measure $P$ on $\Omega$, let $\rho_P^2(f,g) = \mathbb{E}_P(f-g)^2 - (\mathbb{E}_P(f-g))^2$, and for every $\delta > 0$, set $F_\delta = \{f-g\,|\,f,g \in F,\ \rho_P(f,g) \le \delta\}$.

THEOREM 2.11 ([8]).   *$F$ is a uniform Donsker class if and only if the following holds: for every probability measure $P$ on $\Omega$, $G_P$ has a version with bounded, $\rho_P$-uniformly continuous sample paths, and for these versions,*

$$\sup_P \mathbb{E}\sup_{f\in F}|G_P(f)| < \infty, \qquad \lim_{\delta\to 0}\sup_P \mathbb{E}\sup_{h\in F_\delta}|G_P(h)| = 0.$$

The main tool in the analysis of uniform Donsker classes is the Koltchinskii–Pollard entropy integral.

THEOREM 2.12 ([8]).   *If $F \subset B(L_\infty(\Omega))$ satisfies that*

$$\int_0^\infty \sup_n \sup_{\mu_n}\sqrt{\log N(\varepsilon, F, L_2(\mu_n))}\,d\varepsilon < \infty,$$

*then it is a uniform Donsker class.*

**3. Shattering by $B_{X^*}$.**   The goal of this section is to bound the shattering dimension of the dual unit ball of a given Banach space. To that end, we present the geometric interpretation of the shattering dimension when $\Omega \subset X$ and $F = B_{X^*}$.

LEMMA 3.1.   *Let $X$ be a Banach space. Assume that $\{x_1,\dots,x_n\}$ is $\varepsilon$-shattered by $B_{X^*}$ and set $E = \text{span}\{x_1,\dots,x_n\}$. If $A$ is the symmetric convex hull of $\{x_1,\dots,x_n\}$, then $\varepsilon(B_X \cap E) \subset A$.*

PROOF.   Let $\{x_1,\dots,x_n\}$ be $\varepsilon$-shattered by $B_{X^*}$ and let $\{s_1,\dots,s_n\}$ to be a witness to the shattering. Put $(a_i)_{i=1}^{n} \subset \mathbb{R}$, set $I = \{i\,|\,a_i \ge 0\}$ and let $x_I^*$ be the functional shattering the set $I$. For every such $I$ and every $i \in I$,

$$x_I^*(x_i) - x_{I^c}^*(x_i) \ge s_i + \varepsilon - (s_i - \varepsilon) = 2\varepsilon,$$

and if $i \notin I$,

$$x_I^*(x_i) - x_{I^c}^*(x_i) \le s_i - \varepsilon - (s_i + \varepsilon) = -2\varepsilon.$$

Thus,

$$\left\|\sum_{i=1}^{n} a_i x_i\right\| = \sup_{x^* \in B_{X^*}} \left| x^*\left(\sum_{i=1}^{n} a_i x_i\right) \right|$$

$$\geq \tfrac{1}{2} \sup_{x^*, \tilde{x}^* \in B_{X^*}} \left| x^*\left(\sum_{i=1}^{n} a_i x_i\right) - \tilde{x}^*\left(\sum_{i=1}^{n} a_i x_i\right) \right| = (*).$$

Selecting $x^* = x_I^*$ and $\tilde{x}^* = x_{I^c}^*$,

$$(*) \geq \tfrac{1}{2} \left| x_I^*\left(\sum_{i \in I} a_i x_i + \sum_{i \in I^c} a_i x_i\right) - x_{I^c}^*\left(\sum_{i \in I} a_i x_i + \sum_{i \in I^c} a_i x_i\right) \right|$$

$$= \tfrac{1}{2} \left| \sum_{i \in I} a_i (x_I^*(x_i) - x_{I^c}^*(x_i)) + \sum_{i \in I^c} (-a_i)(x_{I^c}^*(x_i) - x_I^*(x_i)) \right|$$

$$\geq \varepsilon \sum_{i=1}^{n} |a_i|.$$

Since every point $x$ on the boundary of $A$ is given by $x = \sum_{i=1}^{n} a_i x_i$, where $\sum_{i=1}^{n} |a_i| = 1$, then $\|x\| = |\sum_{i=1}^{n} a_i x_i| \geq \varepsilon$, which proves our claim. $\square$

COROLLARY 3.2. *The set $\{x_1, \ldots, x_n\} \subset B_X$ is $\varepsilon$-shattered by $B_{X^*}$ if and only if $(x_i)_{i=1}^{n}$ are linearly independent and $\varepsilon$-dominate the $\ell_1^n$ unit-vector basis; that is, for every $a_1, \ldots, a_n \in \mathbb{R}$, $\varepsilon \sum_{i=1}^{n} |a_i| \leq \|\sum_{i=1}^{n} a_i x_i\|$).*

PROOF. Let $E = \operatorname{span}\{x_1, \ldots, x_n\}$ for some linearly independent elements of $B_X$, define $T : \ell_1^n \to \ell_2^n$ by $T e_i = x_i$ and set $A$ to be the symmetric convex hull of $\{x_1, \ldots, x_n\}$. For every $I \subset \{1, \ldots, n\}$, there is some $v \in B_\infty^n$ such that $\langle v, e_i \rangle = 1$ if $i \in I$ and $\langle v, e_j \rangle = -1$ otherwise. Note that $\langle v, e_i \rangle = \langle v, T^{-1} T e_i \rangle = \langle T^{-1*} v, T e_i \rangle$ and that $A^o = (T B_1^n)^o = T^{-1*} B_\infty^n$, implying that $T^{-1*} v \in A^o$. If $\{x_1, \ldots, x_n\}$ $\varepsilon$-dominate the $\ell_1^n$ unit-vector basis, then $\varepsilon(B_X \cap E) \subset A$ and $A^o \subset \varepsilon^{-1}(B_X \cap E)^o = \varepsilon^{-1} P_E B_{X^*}$, where $P_E$ is the orthogonal projection onto $E$. Thus, there is some $x^* \in B_{X^*}$ such that $T^{-1*} v = t P_E x^*$ for some $0 < t \leq \varepsilon^{-1}$. Hence, $\langle x^*, x_i \rangle = \langle x^*, T e_i \rangle = \langle P_E x^*, T e_i \rangle = t^{-1} \langle T^{-1*} v, T e_i \rangle \geq \varepsilon$ if $i \in I$. By a similar argument, $\langle x^*, T e_j \rangle \leq -\varepsilon$ if $j \notin I$, which shows that $\{x_1, \ldots, x_n\}$ is $\varepsilon$-shattered by $B_{X^*}$.

Conversely, if $\{x_1, \ldots, x_n\} \subset B_X$ is $\varepsilon$-shattered, then for every $a_1, \ldots, a_n \in \mathbb{R}$,

$$\varepsilon \sum_{i=1}^{n} |a_i| \leq \left\|\sum_{i=1}^{n} a_i x_i\right\|.$$

Hence, $(x_i)_{i=1}^{n}$ are independent and $\varepsilon$-dominate the $\ell_1^n$ unit-vector basis. $\square$

This result enables us to estimate the shattering dimension of the dual unit ball of an infinite-dimensional Banach space $X$ when considered as a class of functions on $B_X$. It turns out that the shattering dimension is determined by the *type* of $X$.

DEFINITION 3.3.    A Banach space $X$ has type $p$ if there is some constant $C$ such that for every integer $n$ and every $x_1, \ldots, x_n \in X$,

$$(3.1) \qquad \mathbb{E} \left\| \sum_{i=1}^{n} \varepsilon_i x_i \right\| \le C \left( \sum_{i=1}^{n} \|x_i\|^p \right)^{1/p},$$

where $(\varepsilon_i)_{i=1}^{n}$ are independent Rademacher random variables. The smallest constant for which (3.1) holds is called the type $p$ constant of $X$ and is denoted by $T_p(X)$.

The basic facts concerning the concept of type may be found, for example, in [18]. Clearly, for every Banach space (3.1) holds in the case $p = 1$ with $T_1(X) = 1$. If $p^* = \sup\{p | X \text{ has type } p\}$, then $1 \le p^* \le 2$, and if $p^* = 1$, then $X$ is said to have a trivial type.

Recall that the distance between two isomorphic Banach spaces $X$ and $Y$ is defined as $d(X, Y) = \inf \|T\| \cdot \|T^{-1}\|$, where the infimum is taken with respect to all isomorphisms between $X$ and $Y$. It is easy to see that if $X$, $Y$ and $Z$ are isomorphic, then $d(X, Z) \le d(X, Y) \cdot d(Y, Z)$.

THEOREM 3.4.    *Let $X$ be an infinite-dimensional Banach space. Then* $\mathrm{VC}(\varepsilon, B_{X^*}, B_X)$ *is finite for every $\varepsilon > 0$ if and only if $X$ has a nontrivial type. If $X$ has type $p$, then*

$$\left( \frac{1}{\varepsilon} \right)^{p^*/(p^*-1)} - 1 \le \mathrm{VC}(\varepsilon, B_{X^*}, B_X) \le \left( \frac{T_p(X)}{\varepsilon} \right)^{p/(p-1)} + 1.$$

The lower bound and a weaker version of the upper one were established in [13]. We repeat the proof of the lower bound for the sake of completeness.

PROOF OF THEOREM 3.4.    If $\{x_1, \ldots, x_n\}$ is $\varepsilon$-shattered, then it $\varepsilon$-dominates the $\ell_1^n$ unit-vector basis. By selecting $a_i = \varepsilon_i$, $\varepsilon n \le \|\sum_{i=1}^{n} \varepsilon_i x_i\|$. On the other hand, taking the expectation with respect to the Rademacher variables,

$$\mathbb{E} \left\| \sum_{i=1}^{n} \varepsilon_i x_i \right\|_X \le T_p(X) \left( \sum_{i=1}^{n} \|x_i\|_X^p \right)^{1/p} \le T_p(X) n^{1/p}.$$

Thus, there is a realization $(\varepsilon_i)_{i=1}^{n}$ such that $\|\sum_{i=1}^{n} \varepsilon_i x_i\| \le T_p(X) n^{1/p}$. Combining the two inequalities, $n \le (T_p(X)/\varepsilon)^{p/(p-1)}$.

Conversely, for every $\lambda > 0$ and every integer $n$, there is a subspace $X_n \subset X$ such that $\dim X_n = n$ and $d(\ell_{p^*}^n, X_n) \leq 1 + \lambda$ (see [18]). Recall that $d(\ell_1^n, \ell_{p^*}^n) = n^{1-1/p^*}$ (see [21]), hence, $d(X_n, \ell_1^n) \leq (1 + \lambda)n^{1-1/p^*}$, and, in particular, there are $x_1, \ldots, x_n \subset B_X$ such that for every $(a_i)_{i=1}^n \subset \mathbb{R}$,

$$\frac{1}{(1+\lambda)n^{1-1/p^*}} \sum_{i=1}^n |a_i| \leq \left\| \sum_{i=1}^n a_i x_i \right\|.$$

Therefore, $\{x_1, \ldots, x_n\}$ is $n^{(1-p^*)/p^*}(1 + \lambda)^{-1}$-shattered by $B_{X^*}$, and the claim follows by taking $\lambda \to 0$.

The assertion in the case $p^* = 1$ follows in a similar manner. $\quad\square$

The uGC part of the next corollary was first proved in [7], and the second part may also be known to experts; the proof presented below is new, as far as we know.

COROLLARY 3.5. *Let $X$ be an infinite-dimensional Banach space. Then, $F = B_{X^*}$ is a uGC class on $\Omega = B_X$ if and only if $X$ has a nontrivial type. Also, for any infinite-dimensional $X$, $F$ is not a uniform Donsker class on $\Omega$.*

PROOF. The fact that the pair is a uGC class if and only if $X$ has a nontrivial type follows from Theorem 3.4 and the characterization of uGC classes as classes with a finite shattering dimension at every scale $\varepsilon$ (see [1]).

As for the second part, in [8], Example 3.3, it was shown that if $X = \ell_2$, then $F = B_2$ is not a uniform Donsker class on $\Omega = B_2$. Moreover, an easy modification of the proof reveals the following: if there is a constant $C$ such that for every integer $n$ there are spaces $X_n \subset X$ of dimension $n$ for which $d(X_n, \ell_2^n) \leq C$, then $F = B_{X^*}$ is not a uniform Donsker class on $\Omega = B_X$. By Dvoretzky's theorem [18], every infinite-dimensional Banach space has such subspaces $X_n$ (with a constant $C$ arbitrarily close to 1). $\quad\square$

Unlike the infinite-dimensional case, in which the growth of $\mathrm{VC}(\varepsilon, B_{X^*}, B_X)$ is determined by the type of $X$, it is not clear whether the same holds for finite-dimensional spaces; indeed, the lower bound in Theorem 3.4 is based on the fact that $X$ contains spaces which are arbitrarily close to $\ell_{p^*}^n$ for every integer $n$, which is only true for infinite-dimensional spaces.

**4. The shattering dimension of finite-dimensional bodies.** It turns out that some applications require that the set of functionals $F$ is not the dual of the domain but some other convex symmetric set; thus, in the finite-dimensional context it is natural to investigate the following question.

Question 4.1. *Let $K$ and $L$ be two convex symmetric bodies in $\mathbb{R}^d$ and view the elements of $L^\circ$ as functions on $K$ using the fixed inner product in $\mathbb{R}^d$. What is $\mathrm{VC}(\varepsilon, L^\circ, K)$?*

We have shown that $\mathrm{VC}(\varepsilon, L^\circ, K) = n$ if and only if $n$ is the largest such that there are $n$ points $\{x_1, \ldots, x_n\} \subset K$ for which $\varepsilon(L \cap E) \subset \mathrm{absconv}(x_1, \ldots, x_n)$, where $E = \mathrm{span}\{x_1, \ldots, x_n\}$.

The next theorem provides a general upper bound on $\mathrm{VC}(\varepsilon, L^\circ, K)$ based on a volumetric argument. The result is presented for finite-dimensional bodies but can be easily extended to the infinite-dimensional case.

Theorem 4.2. *There is an absolute constant $C$ such that for every two integers $n \leq m$ and every two balls $K, L \subset \mathbb{R}^m$ the following holds: if $\{x_1, ..., x_n\} \subset K$ is $\varepsilon$-shattered by $L^\circ$, then*

$$\sqrt{n} \leq \frac{C}{\varepsilon} \mathrm{vr}((K \cap E)^o) \frac{|K \cap E|^{1/n}}{|L \cap E|^{1/n}},$$

*where $E = \mathrm{span}\{x_1, \ldots, x_n\}$.*

Proof. Assume that $\{x_1, \ldots, x_n\} \subset K$ is $\varepsilon$-shattered by $L^o$. By Lemma 3.1, $\varepsilon(L \cap E) \subset A \subset K \cap E$, where $A$ is the symmetric convex hull of $\{x_1, \ldots, x_n\}$, and, thus, $(K \cap E)^o \subset A^o$. By Lemma 2.8,

$$c\sqrt{n} \leq \mathrm{vr}((K \cap E)^o) \mathrm{cr}((K \cap E)^o) \leq \mathrm{vr}((K \cap E)^o) \left( \frac{|A^o|}{|(K \cap E)^o|} \right)^{1/n}$$

$$\leq \frac{1}{\varepsilon} \mathrm{vr}((K \cap E)^o) \left( \frac{|(L \cap E)^o|}{|(K \cap E)^o|} \right)^{1/n} \leq \frac{C}{\varepsilon} \mathrm{vr}((K \cap E)^o) \left( \frac{|K \cap E|}{|L \cap E|} \right)^{1/n},$$

where the last inequality follows from the Santaló and inverse Santaló inequalities. $\qquad\square$

Combining this theorem with Remark 2.3 on the ratio $|B_p \cap E|/|B_q \cap E|$ and Theorem 2.7 on the volume ratio of projections of $\ell_p$, the following is evident:

Corollary 4.3. *For every $1 \leq p \leq q < \infty$ there is a constant $C_{p,q}$ for which the following holds: if $\{x_1, \ldots, x_n\} \subset B_p$ is $\varepsilon$-shattered by $B_{q'}$, then*

$$\varepsilon \leq C_{p,q} \begin{cases} n^{1/q-1}, & \text{if } 1 \leq p \leq 2, \\ n^{1/q-1/p-1/2}, & \text{if } 2 < p < \infty. \end{cases}$$

In the sequel we will show that this estimate is sharp. Since a similar argument is used in the proof of Theorem 4.10, we shall not present the proof of the corollary here.

Let us mention the following observations: first, using Santaló's inequality, $\mathrm{vr}((K \cap E)^o)|K \cap E|^{1/n} \leq |\tilde{\mathcal{E}}_{K \cap E}|^{1/n}$. Therefore, from the volumetric point of view, all that matters is the ratio between the volume of the ellipsoid of minimal volume containing the section of $K$ spanned by $\{x_1, \ldots, x_n\}$ and the volume of $L \cap E$.

Second, estimating the shattering dimension is equivalent to understanding the behavior of its formal inverse, which, for a given linearly independent set $\{x_1, \ldots, x_n\} \subset K$, is the largest $\varepsilon > 0$ such that $\varepsilon(L \cap E) \subset \mathrm{absconv}(x_1, \ldots, x_n)$, where $E = \mathrm{span}\{x_1, \ldots, x_n\}$. Thus, one can take $K = TB_1^n$, where $T \colon \ell_1^n \to \ell_2$ is defined by $Te_i = x_i$, and the volume of the ellipsoid of minimal volume containing $TB_1^n$ is the significant quantity.

Finally, if $(\lambda_i)_{i=1}^n$ are the singular values of the operator $T$, that is, the eigenvalues of $\sqrt{T^*T}$, then $|\tilde{\mathcal{E}}_{TB_1^n}|^{1/n}$ is equivalent to $n^{-1/2}(\prod_{i=1}^n \lambda_i)^{1/n}$.

4.1. *Shattering and factorization through $\ell_1^n$.* An alternative way to formulate the problem of estimating the shattering dimension is as a factorization problem.

DEFINITION 4.4. For every two balls $K$ and $L$ in $\mathbb{R}^m$ and every integer $n \leq m$, let

$$\Gamma_n(K, L) = \inf \|A\|\|B\|;$$

the infimum is taken with respect to all subspaces of $E \subset \mathbb{R}^m$ of dimension $n$, and all operators $B \colon (E, \|\cdot\|_{L \cap E}) \to \ell_1^n$, $A \colon \ell_1^n \to (E, \|\cdot\|_{K \cap E})$ such that $AB = id \colon L \cap E \to K \cap E$.

The following lemma shows that $1/\Gamma_n(K, L)$ is the formal inverse of the shattering dimension.

LEMMA 4.5. *For every integer $n$ and any balls $K$ and $L$,*

$$\text{(4.1)} \quad \frac{1}{\Gamma_n(K, L)} = \sup\{\varepsilon | \exists \{x_1, \ldots, x_n\} \subset K, \ \varepsilon(L \cap E) \subset \mathrm{absconv}(x_1, \ldots, x_n)\}$$

$$\text{(4.2)} \quad = \sup\{\varepsilon | \ \mathrm{VC}(\varepsilon, L^\circ, K) \geq n\},$$

*where $E = \mathrm{span}\{x_1, \ldots, x_n\}$.*

PROOF. If the identity admits an optimal factorization $id = AB$, set $A' = A/\|A\|_{\ell_1^n \to K \cap E}$ and observe that the set $\{A'e_1, \ldots, A'e_n\} \subset K \cap E$ sat-isfies that for any $a_1, \ldots, a_n \in \mathbb{R}$,

$$\|A\|_{\ell_1^n \to K \cap E} \cdot \|B\|_{L \to E \to \ell_1^n} \cdot \left\|\sum_{i=1}^n a_i A' e_i\right\|_L \geq \left\|B\left(\sum_{i=1}^n a_i A e_i\right)\right\|_{\ell_1^n}$$

$$\geq \left\|\sum_{i=1}^n a_i e_i\right\|_{\ell_1^n} = \sum_{i=1}^n |a_i|.$$

Hence, absconv$(A' e_1, \ldots, A' e_n) \subset K \cap E$ contains $(\|A\|\|B\|)^{-1}(L \cap E)$ and

$$\frac{1}{\Gamma_n(K, L)} \leq \sup\{\varepsilon | \exists \{x_1, \ldots, x_n\} \subset K, \ \varepsilon(L \cap E) \subset \text{absconv}(x_1, \ldots, x_n)\}.$$

For the reverse inequality, if $\{x_1, \ldots, x_n\} \subset K$ are such that $\varepsilon(L \cap E) \subset$ absconv$(x_1, \ldots, x_n)$, define $T \colon \mathbb{R}^n \to \mathbb{R}^n$ by $T e_i = x_i$. Clearly, $\|T\|_{\ell_1^n \to K \cap E} \leq 1$ and

$$\|T^{-1}\|_{L \cap E \to \ell_1^n} = \sup_{x \in L \cap E} \|T^{-1} x\|_{\ell_1^n} = \sup_{x \in L \cap E} \|x\|_{T \ell_1^n} \leq \frac{1}{\varepsilon}.$$

Thus, $\|T\|_{\ell_1^n \to K \cap E} \cdot \|T^{-1}\|_{L \cap E \to \ell_1^n} \leq 1/\varepsilon$ and $1/\Gamma_n(K, L) \geq \varepsilon$. $\quad\square$

Combining Theorem 4.2 and Lemma 4.5, we obtain the following:

COROLLARY 4.6. *There is an absolute constant $c > 0$ such that for any two integers $n \leq m$ and any two balls $K, L \subset \mathbb{R}^m$,*

$$\Gamma_n(K, L) \geq c \frac{\sqrt{n}}{v_n(id \colon K \to L) \sup_E \text{vr}(P_E K^\circ)},$$

*where $\dim(E) = n$.*

4.2. *Factorization constants of $\ell_p^m$.* The goal of the next section is to investigate the shattering dimension of the class of linear functionals $F = B_{q'}^m$ on $\Omega = B_p^m$ for $1 \leq p, q \leq \infty$. First, in Theorems 4.9 and 4.10 below we present a tight estimate on the factorization constant of $id \colon \ell_q^n \to \ell_p^n$ through $\ell_1^n$. Then, we use this result to estimate $\Gamma_n(B_p^m, B_q^m)$ and, thus, bound VC$(\varepsilon, B_{q'}^m, B_p^m)$; finally, we show that if $1 \leq p < q \leq \infty$, then $F = B_{q'}$ is a uniform Donsker class on $\Omega = B_p$.

We begin with two lemmas needed for the proof of Theorem 4.9.

LEMMA 4.7. *Let $\mu$ be the Haar measure on the $n$ dimensional sphere $S^{n-1}$. Set $K$ and $L$ to be balls in $\mathbb{R}^n$ and put $\alpha$ to be such that*

$$\mu\left(x \in S^{n-1} \middle| \|x\|_K > \frac{1}{\alpha}\right) < \frac{1}{2n}.$$

*If $\varepsilon$ satisfies that*

$$\mu\left(x \in S^{n-1} \middle| \|x\|_{L^\circ} > \frac{\alpha}{\varepsilon\sqrt{n}}\right) < 2^{-(n+1)},$$

*then $\Gamma_n(K, L) \leq 1/\varepsilon$.*

PROOF. Denote by $O_n$ the orthogonal group and let $P_{O_n}$ be the Haar measure on $O_n$. Set $U \in O_n$ and define $x_i = \alpha U e_i$. Using the standard connection between $P_{O_n}$ and $\mu$ on $S^{n-1}$,

$$P_{O_n}(x_i \in K) = \mu\left(x \in S^{n-1} \middle| \|x\|_K \leq \frac{1}{\alpha}\right),$$

hence,

$$P_{O_n}(x_i \in K \text{ for all } i) \geq 1 - n\mu\left(x \in S^{n-1} \middle| \|x\|_K > \frac{1}{\alpha}\right) > \frac{1}{2}.$$

Moreover,

$$P_{O_n}(\text{conv}(\pm \alpha U e_i) \supset \varepsilon L) = P_{O_n}\left(\sup_{(\sigma_i)_{i=1}^n \in \{-1,1\}^n} \left\|\frac{1}{\alpha}\sum_{i=1}^n \sigma_i U e_i\right\|_{L^\circ} \leq \frac{1}{\varepsilon}\right).$$

For every vector $(\sigma_1, \ldots, \sigma_n) \in \{-1, 1\}^n$,

$$P_{O_n}\left(\left\|\frac{1}{\alpha}\sum_{i=1}^n \sigma_i U e_i\right\|_{L^\circ} > \frac{1}{\varepsilon}\right) = P_{O_n}\left(\left\|U\left(\frac{1}{\sqrt{n}}\sum_{i=1}^n e_i\right)\right\|_{L^\circ} > \frac{\alpha}{\varepsilon\sqrt{n}}\right)$$

$$= \mu\left(x \in S^{n-1} \middle| \|x\|_{L^\circ} > \frac{\alpha}{\varepsilon\sqrt{n}}\right).$$

Thus,

$$P_{O_n}(\text{conv}(\pm \alpha U e_i) \supset \varepsilon L) \geq 1 - 2^n \mu\left(x \in S^{n-1} \middle| \|x\|_{L^\circ} \geq \frac{\alpha}{\varepsilon\sqrt{n}}\right) > \frac{1}{2},$$

and there is some orthogonal operator which belongs to both events. The operator $T = \alpha U$ satisfies that $\|T\|_{\ell_1^n \to K \cap E} \leq 1$ and $\|T^{-1}\|_{L \cap E \to \ell_1^n} \leq 1/\varepsilon$, as claimed. $\square$

LEMMA 4.8. *There are constants $C_p$ for which the following holds: for every integer $n$,*

$$\mu(x \in S^{n-1} | \|x\|_{\ell_p^n} \geq C_p n^{1/p - 1/2}) \leq 2^{-(n+1)},$$

*if $1 \leq p \leq 2$ and if $2 \leq p < \infty$, then*

$$\mu(x \in S^{n-1} | \|x\|_{\ell_p^n} \geq C_p n^{1/p - 1/2}) \leq e^{-n^{2/p}}.$$

PROOF.   Denote by $M(B_p^n)$ the median of $\|x\|_p$ on $S^{n-1}$. By Lévy's inequality [18],

$$\mu(x \in S^{n-1} | \|x\|_p \geq (1+t)M(B_p^n)) \leq \exp\left\{-\frac{t^2 n M^2(B_p^n)}{2\|id\|_{\ell_2^n \to \ell_p^n}^2}\right\}.$$

Recall that $M(B_p^n) \sim_p n^{1/p-1/2}$ (see, e.g., [18]) and that $\|id\|_{\ell_2^n \to \ell_p^n} = \max\{n^{1/p-1/2}, 1\}$. It follows that for $1 \leq p \leq 2$ and $C$ large enough, depending only on $p$,

$$\mu(x \in S^{n-1} | \|x\|_{\ell_p^n} \geq C n^{1/p-1/2}) \leq e^{-c_p C^2 n} \leq 2^{-(n+1)},$$

while for $2 \leq p < \infty$ and $C$ depending only on $p$,

$$\mu(x \in S^{n-1} | \|x\|_{\ell_p^n} \geq C n^{1/2-1/p}) \leq e^{-n^{2/p}}. \qquad \square$$

The above results will play an important role in the proof of the following theorem, in which we construct factorizations of $id: \ell_q^n \to \ell_p^n$ through $\ell_1^n$.

THEOREM 4.9.   Let $K = B_p^n$ and $L = B_q^n$. Then, $\Gamma_n(K,L)$ satisfies that

$$\Gamma_n(K,L) \leq C_{p,q} \begin{cases} n^{1/2+1/p-1/q}, & \text{if } 2 \leq p,q \leq \infty, \\ n^{1-1/q}, & \text{if } 1 \leq p \leq 2, \\ n^{1-1/q}, & \text{if } 1 \leq q \leq 2 \leq p \leq \infty \text{ and } p' > q, \\ n^{1/p}, & \text{if } 1 \leq q \leq 2 \leq p \leq \infty \text{ and } p' \leq q. \end{cases}$$

PROOF.   First, assume that $2 \leq p < \infty$ and $2 \leq q \leq \infty$. By Lemmas 4.8 and 4.7 it suffices to choose $1/\alpha = C_p n^{1/p-1/2}$ and select $\varepsilon$ which satisfies that $C_q n^{1/q'-1/2} = C_q n^{1/2-1/q} = C_q \alpha/\varepsilon\sqrt{n}$, that is, $\frac{1}{\varepsilon} \sim_q \frac{n^{1-1/q}}{\alpha}$. Therefore, $\Gamma_n(B_p^n, B_q^n) \leq C_{p,q} n^{1/2+1/p-1/q}$.

Next, if $1 \leq p \leq 2$, then

$$\Gamma_n(K,L) \leq \|id\|_{\ell_q^n \to \ell_1^n} \|id\|_{\ell_1^n \to \ell_p^n} = n^{1-1/q}.$$

If $2 \leq p < \infty$ and $1 \leq q \leq 2$, one has to treat two cases; if $1 \leq q \leq p'$, then using the identity operator as above, $\Gamma_n(B_p^n, B_q^n) \leq n^{1-1/q}$. On the other hand, if $p' \leq q \leq 2$, then by the first part of our claim,

$$\Gamma_n(B_p^n, B_q^n) \leq \|id\|_{\ell_q^n \to \ell_2^n} \Gamma_n(B_p^n, B_2^n) \leq C_p n^{1/p}.$$

Finally, one has to address the situation when $p$ is infinity. If $p = q = \infty$, then $\Gamma_n(B_\infty^n, B_\infty^n) = d(\ell_1^n, \ell_\infty^n) \leq C n^{1/2}$ [21].

For $p = \infty$ we first examine the case $q = 2$. Let $\ell_p^n(\mathbb{C})$ to be $\mathbb{C}^n$ endowed with the $\ell_p$ norm and set $T = (n^{-1/2} e^{2\pi ijk/n})_{j,k=1}^n$. It is easy to check that $\|T\|_{\ell_1^n(\mathbb{C}) \to \ell_\infty^n(\mathbb{C})} \leq n^{-1/2}$ and that $\|T\|_{\ell_\infty^n(\mathbb{C}) \to \ell_2^n(\mathbb{C})} \leq n^{1/2}$. For our purpose, $\ell_p^n(\mathbb{C})$ can be considered as the $\ell_p^n$ sum of two-dimensional Euclidean spaces, $\ell_2^2$, over the reals, and since for any $1 \le p \le \infty$, $\|id\|_{\ell_p^n(\mathbb{C}) \to \ell_p^{2n}} \cdot \|id\|_{\ell_p^{2n} \to \ell_p^n(\mathbb{C})} \le \sqrt{2}$, then $\Gamma_{2n}(\ell_2^{2n}, \ell_\infty^{2n}) \le 2$. The case where $n$ is odd is easily reduced to the even case.

Finally, for a general $q$,

$$\Gamma_n(B_\infty^n, B_q^n) \le \|id\|_{\ell_q^n \to \ell_2^n} \Gamma_n(B_\infty^n, B_2^n) \le C \|id\|_{\ell_q^n \to \ell_2^n},$$

as claimed. $\square$

The next step in our analysis is to show that the bounds in Theorem 4.9 are tight. The proof uses the notion of $r$-summing operators. Recall that an operator $T \colon X \to Y$ is $r$-summing for $1 \le r < \infty$, if there is a $C < \infty$ such that

$$(4.3) \qquad \sum_{i=1}^n \|T x_i\|^r \le C^r \sup_{x^* \in B_{X^*}} \sum_{i=1}^n |x^*(x_i)|^r$$

for all integers $n$ and all $x_1, \dots, x_n \in X$. The smallest $C$ for which (4.3) holds is denoted by $\pi_r(T)$ and is called the $r$-summing norm of $T$.

**Theorem 4.10.** *There exist $c_{p,q}$ such that if $K = B_p^n$ and $L = B_q^n$, then $\Gamma_n(K, L)$ satisfies that*

$$\Gamma_n(K, L) \ge c_{p,q} \begin{cases} n^{1/2 + 1/p - 1/q}, & \text{if } 2 \le p, q \le \infty, \\ n^{1 - 1/q}, & \text{if } 1 \le p \le 2. \end{cases}$$

*Also, there are $c_{p,q,r}$ such that if $1 < q \le 2 \le p < \infty$ and $r > \max\{p, q'\}$, then*

$$\Gamma_n(K, L) \ge c_{p,q,r} n^{1/r}.$$

**Proof.** The first two cases follow from the volumetric estimate as in Corollary 4.3. Indeed, if $\{x_1, \dots, x_n\} \subset K$ is $\varepsilon$-shattered by $L^\circ$, then

$$\varepsilon \sqrt{n} \le C \operatorname{vr}(K^\circ) \left( \frac{|K|}{|L|} \right)^{1/n}$$

for some absolute constant $C$. Since $K = B_p^n$ and $L = B_q^n$, then $(|K|/|L|)^{1/n} \sim_{p,q} n^{1/q - 1/p}$ and $\operatorname{vr}(K^\circ) \sim n^{1/p - 1/2}$ for $1 \le p \le 2$ and $\sim 1$ for $2 \le p \le \infty$. Hence,

$$\varepsilon \le C_{p,q} \begin{cases} n^{1/q - 1/p - 1/2}, & \text{if } 2 \le p \le \infty, \\ n^{1/q - 1}, & \text{if } 1 \le p \le 2, \end{cases}$$

and the lower estimate on $\Gamma_n$ is evident from Lemma 4.5.

For $1 \le q < 2 < p \le \infty$ we can get a better estimate than what the volumetric estimates provide. We first investigate $id \colon \ell_q^n \to \ell_{q'}^n$. Observe that if $AB$ is a factorization of the identity through $\ell_1^n$, then $B^*A^*$ is a factorization of $id$ through $\ell_\infty^n$. A theorem of Maurey (see [21], Theorem 21.4(ii)) asserts that, for every $r > q'$, $B^*$ is $r$-summing with $\pi_r(B^*) \le C_{q,r}\|B^*\|$ and, thus, by the properties of $\pi_r$, $\pi_r(id) \le \|A^*\|\pi_r(B^*) \le C_{q,r}\Gamma_n(B_{q'}, B_q)$.

The behavior of the $\pi_r$ norm of the identity between $\ell_p^n$ and other spaces was investigated in [9] and [10]. In particular, in the range we are interested in, it is proved in [9] that $\pi_r(id \colon \ell_q^n \to \ell_{q'}^n) \ge c_{q,r}n^{1/r}$. (For the interested reader, we found that the best way to understand this is to apply Theorem 1 there to our setup. This is rather easy, as is the proof of Theorem 1.)

This settles the case $p = q'$. Turning to the general case, assume first that $2 \le q' \le p < r < \infty$. For any factorization $AB = id_{q \to p}$, $id_{p \to q'}AB$ is a factorization of $id_{q \to q'}$. Therefore, for any $s > q'$,

$$C_{q,r}n^{1/s} \le \|B\|\|id_{p \to q'}A\| \le \|A\|\|B\|\|id_{p \to q'}\| \le \|A\|\|B\|n^{1/q'-1/p},$$

hence,

$$\|A\|\|B\| \ge C_{q,r}n^{1/p+1/s-1/q'}.$$

Choosing $s$ such that $1/r = 1/p + 1/s - 1/q'$ gives the result in this case.

A similar argument may be used to handle the case $q' > p$.  $\square$

Next we estimate $\Gamma_n(B_p^m, B_q^m)$ when $n \le m$. Note that the results we obtain are not for the full range of $p$ and $q$.

THEOREM 4.11.   *For every integers $n \le m$ the following holds:*

1. *If $2 \le q \le p < \infty$ then $\Gamma_n(B_p^m, B_q^m) \sim_{p,q} n^{1/2}m^{1/p-1/q}$.*
2. *If $q \le p \le 2$ then $\Gamma_n(B_p^m, B_q^m) \sim_{p,q} n^{1-1/p}m^{1/p-1/q}$.*
3. *If $p \le q$ and $1 \le p \le 2$ then $\Gamma_n(B_p^m, B_q^m) \sim_{p,q} n^{1-1/q}$.*
4. *If $p \le q$ and $2 < p < \infty$ then $\Gamma_n(B_p^m, B_q^m) \sim_{p,q} n^{1/2+1/p-1/q}$.*

PROOF.   In all cases, the lower bound follows from Corollary 4.6 combined with the estimate on the volume numbers of $id_{p \to q}$ in (2.2) and (2.4), and the volume ratios of quotients of $\ell_p^m$ from Theorem 2.7.

As for the upper bound, the optimal choice in (1) and (2) (at least when $n$ divides $m$) is the section $E$ spanned by

$$v_j = \sum_{i=1}^k e_{j+k(i-1)}, \qquad j = 1, \ldots, n.$$

Then $B_p^m \cap E = (m/n)^{1/2-1/p}B_p^n$. Clearly, $\Gamma_n(B_p^m, B_q^m) \le \Gamma_n(B_p^m \cap E, B_q^m \cap E)$ and when $2 \le q \le p$ the latter can be approximated using the probabilistic argument from Lemmas 4.8 and 4.7. Indeed, a straightforward computation shows that one can take $\alpha = m^{1/2-1/p}$ and that $\varepsilon$ needs to satisfy that $m^{1/2-1/q} = \alpha/n^{1/2}\varepsilon = m^{1/2-1/p}/n^{1/2}\varepsilon$. Thus, $1/\varepsilon \le n^{1/2}m^{1/p-1/q}$, which proves the bound is tight.

When $q \leq p \leq 2$ one uses the identity operator as the factorizing operator between $(m/n)^{1/2-1/q} B_q^n$ and $(m/n)^{1/2-1/p} B_p^n$ to obtain the required result.

The upper bound in (3) is obtained by taking the canonical section $\mathrm{span}\{e_1, \ldots, e_n\}$ and applying Theorem 4.9. $\square$

Some of the information one can obtain from these estimates is summarized in the following:

THEOREM 4.12. *Let $1 \leq p < q \leq \infty$, set $F = B_{q'}$ and $\Omega = B_p$. Then:*

1.
$$\mathrm{VC}(\varepsilon, F, \Omega) \overset{p,q}{\sim} \begin{cases} \varepsilon^{-q/(q-1)}, & \text{if } 1 \leq p \leq 2, \\ \varepsilon^{-1/(1/2+1/p-1/q)}, & \text{if } 2 < p \leq \infty. \end{cases}$$

2. *$F$ is a uniform Donsker class on $\Omega$.*
3. *There are constants $C_{p,q}$ such that for any probability measure $\mu$ on $B_p$, every integer $n$ and every $t > 0$,*

$$Pr\left\{ \sup_{x^* \in B_{q'}} \left| \mathbb{E}_\mu x^* - \frac{1}{n} \sum_{i=1}^n x^*(X_i) \right| \geq C_{p,q} \left( \frac{1}{\sqrt{n}} + \frac{t}{n} \right) \right\} \leq e^{-t},$$

*where $(X_i)_{i=1}^n$ are independent and distributed according to $\mu$.*

Before presenting the proof, we require an additional lemma which follows from Theorem 1.5. Although the first equality is not needed in the sequel, it might be useful in other applications.

LEMMA 4.13. *For any $1 \leq p < q \leq \infty$ there is a constant $C_{p,q}$ for which the following holds: if $x_1, \ldots, x_n \in B_p$ and $T : \ell_{q'} \to \ell_2^n$ is given by $Tx^* = n^{-1/2} \sum_{i=1}^n x^*(x_i) e_i$, then, for every $\varepsilon > 0$,*

$$\log N(\varepsilon, T B_{q'}, \ell_2^n) = \log N(\varepsilon, B_{q'}, L_2(\mu_n)) \leq C_{p,q} \mathrm{VC}(\varepsilon, B_{q'}, B_p) \cdot \log \frac{2}{\varepsilon},$$

*where $\mu_n$ is the empirical measure supported on $\{x_1, \ldots, x_n\}$.*

PROOF OF THEOREM 4.12. The first part of the claim follows from Corollary 4.3 which yields the upper bound, while the lower one follows immediately from Theorem 4.11.

The second part is evident because, by Lemma 4.13, the class has a converging entropy integral, which by Theorem 2.12 suffices to ensure that $F$ is a uniform Donsker class.

Finally, the last part follows from the first, combined with Talagrand's inequality (Theorem 1.3) and the estimate on the expected deviation in terms of the shattering dimension (Theorem 1.6). $\square$

**5. The shattering dimension of images of $B_1^m$.** Although the volumetric approach yields sharp results in some cases, and, in particular, for $\Gamma_n(B_p^n, B_q^n)$ for a certain range of $p$ and $q$, an exact estimate on the factorization constant $\Gamma_n(TB_1^n, B_q^n)$ does not follow from the volumetric argument, since the position of $B_1^n$ is significant, and not only the volume of the ellipsoid of minimal volume containing $TB_1^n$. Indeed, we show that spectral information does not suffice for sharp estimates on the shattering dimension. To demonstrate this, given a set of (nonnegative) singular values (arranged in a nonincreasing order) $\Lambda = (\lambda_1, \dots, \lambda_n)$, let $\mathbb{T}_\Lambda$ be the subset of $GL_n$ consisting of the matrices which have $\Lambda$ as singular values.

THEOREM 5.1.  *For every set $\Lambda$ of singular values,*

$$\sup_{T \in \mathbb{T}_\Lambda} \Gamma_n(TB_1^n, B_q^n) = \frac{1}{\lambda_n} \begin{cases} n^{1-1/q}, & \text{if } q \geq 2, \\ n^{1/2}, & \text{if } q < 2, \end{cases}$$

*and*

$$\inf_{T \in \mathbb{T}_\Lambda} \Gamma_n(TB_1^n, B_q^n) \stackrel{q}{\sim} \left( \sum_{i=1}^n \lambda_i^{-2} \right)^{1/2} \begin{cases} 1, & \text{if } q \geq 2, \\ n^{1/2-1/q}, & \text{if } q < 2. \end{cases}$$

To compare this result to the one obtained via the volumetric approach (Theorem 4.2), take $q = 2$, and recall that Theorem 4.2 implies that

$$\Gamma_n(TB_1^n, B_2^n) \geq cn^{1/2} \left( \prod_{i=1}^n \lambda_i^{-2} \right)^{1/2n},$$

which, by the means inequality, is weaker than the conclusion of Theorem 5.1.

PROOF OT THEOREM 5.1.  By Lemma 3.1, for every $T \in GL_n$,

$$\Gamma_n(TB_1^n, B_q^n) = (\sup\{\varepsilon | \varepsilon B_q^n \subset TB_1^n\})^{-1} = \max_{\|x\|_q=1} \|x\|_{TB_1^n}.$$

Since $(TB_1^n)^o = T^{-1*}B_\infty^n$, then for every $x$,

$$\|x\|_{TB_1^n} = \sup_{y \in (TB_1^n)^o} \langle x, y \rangle$$

$$= \sup_{y \in B_\infty^n} \langle x, T^{-1*}y \rangle$$

$$= \sup_{(\varepsilon_i)_{i=1}^n \in \{-1,1\}^n} \left\langle T^{-1}x, \sum_{i=1}^n \varepsilon_i e_i \right\rangle$$

and

$$\max_{\|x\|_q=1} \|x\|_{TB_1^n} = \sup_{(\varepsilon_i)_{i=1}^n \in \{-1,1\}^n} \sup_{\|x\|_q=1} \left\langle T^{-1}x, \sum_{i=1}^n \varepsilon_i e_i \right\rangle.$$

By the polar decomposition, $T^{-1} = ODU$, where $V$ and $O$ are orthogonal and $D$ is the diagonal matrix with eigenvalues $\lambda_i^{-1}$. Thus

$$\inf_{T \in \mathbb{T}_\Lambda} \max_{\|x\|_q=1} \|x\|_{TB_1^n}$$

$$= \inf_{O,V \in O_n} \sup_{(\varepsilon_i)_{i=1}^n \in \{-1,1\}^n} \sup_{\|x\|_q=1} \left\langle ODVx, \sum_{i=1}^n \varepsilon_i e_i \right\rangle$$

and

$$\sup_{T \in \mathbb{T}_\Lambda} \max_{\|x\|_q=1} \|x\|_{TB_1^n}$$

$$= \sup_{O,V \in O_n} \sup_{(\varepsilon_i)_{i=1}^n \in \{-1,1\}^n} \sup_{\|x\|_q=1} \left\langle ODVx, \sum_{i=1}^n \varepsilon_i e_i \right\rangle,$$

where $O_n$ denotes the set of orthogonal matrices on $\mathbb{R}^n$. Set $(\mu_i)_{i=1}^n$ to be the eigenvalues of $D$ arranged in a nonincreasing order, that is, $\mu_1 = \lambda_n^{-1} \geq \cdots \geq \mu_n = \lambda_1^{-1}$.

Let

$$f(O,V) = \max_{(\varepsilon_i)_{i=1}^n \in \{-1,1\}^n} \max_{\|x\|_q=1} \left\langle ODVx, \sum_{i=1}^n \varepsilon_i e_i \right\rangle,$$

and observe that

$$f(O,V) = \max_{\|x\|_q=1} \max_{(\varepsilon_i)_{i=1}^n} \left\langle x, \sum_{k=1}^n \left( \sum_{j=1}^n V_{jk} \mu_j \sum_{i=1}^n \varepsilon_i O_{ij} \right) e_k \right\rangle$$

$$= \max_{(\varepsilon_i)_{i=1}^n} \left( \sum_{k=1}^n \left| \sum_{j=1}^n \mu_j \left( \sum_{i=1}^n \varepsilon_i O_{ij} \right) V_{jk} \right|^{q'} \right)^{1/q'}.$$

Clearly,

$$\max_{O,V} f(O,V) = \max_{O,V} \max_{(\varepsilon_i)_{i=1}^n} \left( \sum_{k=1}^n \left| \sum_{j=1}^n \mu_j \left( \sum_{i=1}^n \varepsilon_i O_{ij} \right) V_{jk} \right|^{q'} \right)^{1/q'}$$

$$= \max_V \max_{\|z\|_2=\sqrt{n}} \left( \sum_{k=1}^n \left| \sum_{j=1}^n \mu_j z_j V_{jk} \right|^{q'} \right)^{1/q'}$$

$$= \max_V \max_{\|x\|_2 = \mu_1 \sqrt{n}} \|xV\|_{q'}$$

$$= \mu_1 \sqrt{n} \max_{x \neq 0} \frac{\|x\|_{q'}}{\|x\|_2},$$

from which the first part of the claim follows.

To prove the second part, note that

$$\min_{O,V} (f(O,V))^{q'} \geq \min_{O,V} \mathbb{E}_\varepsilon \sum_{k=1}^n \left| \sum_{j=1}^n \mu_j \left( \sum_{i=1}^n \varepsilon_i O_{ij} \right) V_{jk} \right|^{q'}$$

$$= \min_{O,V} \sum_{k=1}^n \mathbb{E}_\varepsilon \left| \sum_{i=1}^n \varepsilon_i \sum_{j=1}^n \mu_j V_{jk} O_{ij} \right|^{q'} = (*),$$

where $(\varepsilon_i)_{i=1}^n$ are independent Rademacher random variables. Therefore, by Khintchine's inequality,

$$(*) \geq \min_{O,V} C_q \sum_{k=1}^n \left( \sum_{i=1}^n \left( \sum_{j=1}^n \mu_j V_{jk} O_{ij} \right)^2 \right)^{q'/2}.$$

Denoting $h_k = (\mu_j V_{jk})_{j=1}^n$,

$$\left( \sum_{i=1}^n \left( \sum_{j=1}^n \mu_j V_{jk} O_{ij} \right)^2 \right)^{1/2} = \|h_k O\|_2 = \|h_k\|_2 = \left( \sum_{j=1}^n \mu_j^2 V_{jk}^2 \right)^{1/2},$$

and applying Khintchine's inequality again,

$$\left( \sum_{j=1}^n \mu_j^2 V_{jk}^2 \right)^{q'/2} \geq C_q \mathbb{E}_\varepsilon \left| \sum_{j=1}^n \varepsilon_j \mu_j V_{jk} \right|^{q'/2}.$$

By Jensen's inequality and since the matrix $(\varepsilon_j V_{jk})_{j,k=1}^n$ is also orthogonal for any realization of the Rademacher variables,

$$\min_{O,V} f(O,V) \geq C_q \min_V \left( \mathbb{E}_\varepsilon \sum_{k=1}^n \left| \sum_{j=1}^n \varepsilon_j \mu_j V_{jk} \right|^{q'} \right)^{1/q'}$$

$$\geq C_q \mathbb{E}_\varepsilon \min_V \left( \sum_{k=1}^n \left| \sum_{j=1}^n \varepsilon_j \mu_j V_{jk} \right|^{q'} \right)^{1/q'}$$

$$= C_q \min_V \|\mu V\|_{q'}$$

and

$$\min_V \|\mu V\|_{q'} = \|\mu\|_2 \begin{cases} 1, & \text{if } q' \leq 2, \\ n^{1/q' - 1/2}, & \text{if } q' > 2. \end{cases}$$

Finally, to see that the lower bound is tight, set $O = id$, and, thus,

$$f(id, V) = \max_{(\varepsilon_i)_{i=1}^n} \left\| \sum_{i=1}^n \varepsilon_i (V^* \mu)_i e_i \right\|_{q'} = \|V^* \mu\|_{q'}.$$

The sharpness is evident by optimizing with respect to $V$. $\square$

## REFERENCES

[1] Alon, N., Ben-David, S., Cesa-Bianchi, N. and Haussler, D. (1997). Scale sensitive dimensions, uniform convergence and learnability. *J. Assoc. Comput. Mach.* **44** 615–631. MR1481318

[2] Anthony, M. and Bartlett, P. L. (1999). *Neural Network Learning, Theoretical Foundations.* Cambridge Univ. Press. MR1741038

[3] Ball, K. (1989). Volumes of sections of cubes and related problems. *Geometric Aspects of Functional Analysis. Lecture Notes in Math.* **1376** 251–260. Springer, Berlin. MR1008726

[4] Bartlett, P. L., Bousquet, O. and Mendelson, S. (2004). Localized Rademacher complexities. *Ann. Statist.* To appear. MR2040404

[5] Bousquet, O. (2002). A Bennett concentration inequality and its application to suprema of empirical processes. *C. R. Acad. Sci. Paris Ser. I* **334** 495–500. MR1890640

[6] Dudley, R. M. (1999). *Uniform Central Limit Theorems.* Cambridge Univ. Press. MR1720712

[7] Dudley, R. M., Giné, E. and Zinn, J. (1991). Uniform and universal Glivenko–Cantelli classes. *J. Theoret. Probab.* **4** 485–510. MR1115159

[8] Giné, E. and Zinn, J. (1991). Gaussian characterization of uniform Donsker classes of functions. *Ann. Probab.* **19** 758–782. MR1106285

[9] Gluskin, E. D. (1978). Estimates of the norms of certain $p$-absolutely summing operators. *Funktsional. Anal. i Prilozhen.* **12** 24–31. [Translated in *Functional Anal. Appl.* **12** 94–101.] MR498895

[10] Gluskin, E. D., Pietsch, A. and Puhl, J. (1980). A generalization of Khintchine's inequality and its application in the theory of operator ideals. *Studia Math.* **67** 149–155. MR583295

[11] Gordon, Y. and Junge, M. (1997). Volume formulas in $L_p$ spaces. *Positivity* **1** 7–43. MR1659607

[12] Gordon, Y. and Junge, M. (1999). Volume ratios in $L_p$ spaces. *Studia Math.* **136** 147–182. MR1716174

[13] Mendelson, S. (2002). Learnability in Hilbert spaces with reproducing kernels. *J. Complexity* **18** 152–170. MR1895081

[14] Mendelson, S. (2002). Rademacher averages and phase transitions in Glivenko–Cantelli class. *IEEE Trans. Inform. Theory* **48** 251–263. MR1872178

[15] Mendelson, S. (2003). A few notes on statistical learning theory. *Advanced Lectures in Machine Learning. Lecture Notes in Comput. Sci.* **2600** 1–40. Springer, New York.

[16] Mendelson, S. and Vershynin, R. (2003). Entropy and the combinatorial dimension. *Invent. Math.* **152** 37–55. MR1965359

[17] Meyer, M. and Pajor, A. (1988). Sections of the unit ball of $\ell_p^n$. *J. Funct. Anal.* **80** 109–123. MR960226

[18] MILMAN, V. D. and SCHECHTMAN, G. (1986). *Asymptotic Theory of Finite Dimensional Normed Spaces. Lecture Notes in Math.* **1200**. Springer, Berlin. MR856576

[19] PISIER, G. (1989). *The Volume of Convex Bodies and Banach Space Geometry.* Cambridge Univ. Press. MR1036275

[20] TALAGRAND, M. (1994). Sharper bounds for Gaussian and empirical processes. *Ann. Probab.* **22** 28–76. MR1258865

[21] TOMCZAK-JAEGERMANN, N. (1989). *Banach–Mazur Distance and Finite-Dimensional Operator Ideals.* CRC Press, Boca Raton, FL. MR993774

[22] VAN DER VAART, A. W. and WELLNER, J. A. (1996). *Weak Convergence and Empirical Processes.* Springer, New York. MR1385671

[23] VAPNIK, V. and CHERVONENKIS, A. (1971). Necessary and sufficient conditions for uniform convergence of means to mathematical expectations. *Theory Probab. Appl.* **26** 532–553. MR627861

RESEARCH SCHOOL OF INFORMATION
    SCIENCES AND ENGINEERING
AUSTRALIAN NATIONAL UNIVERSITY
CANBERRA ACT 0200
AUSTRALIA
E-MAIL: shahar.mendelson@anu.edu.au

DEPARTMENT OF MATHEMATICS
WEIZMANN INSTITUTE
REHOVOT
ISRAEL
E-MAIL: gideon@weizmann.ac.il